\documentclass[12pt]{article} % use larger type; default would be 10pt
\usepackage{setspace}
\usepackage[utf8]{inputenc} % set input encoding (not needed with XeLaTeX)

\usepackage{subfigure}
\usepackage[margin=1in]{geometry}
\usepackage{graphicx} % support the \includegraphics command and options
\usepackage{epstopdf}

\usepackage{booktabs} % for much better looking tables
\usepackage{array} % for better arrays (eg matrices) in maths
\usepackage{paralist} % very flexible & customisable lists (eg. enumerate/itemize, etc.)
\usepackage{verbatim} % adds environment for commenting out blocks of text & for better verbatim
\usepackage{subfig} % make it possible to include more than one captioned figure/table in a single float
\usepackage{fancyhdr} % This should be set AFTER setting up the page geometry
\usepackage{sectsty}
\allsectionsfont{\sffamily\mdseries\upshape} % (See the fntguide.pdf for font help)
\usepackage[nottoc,notlof,notlot]{tocbibind} % Put the bibliography in the ToC
\usepackage[titles,subfigure]{tocloft} % Alter the style of the Table of Contents

\usepackage{amsmath, amsthm, amssymb}
\usepackage{multirow}
\usepackage{algorithm}
\usepackage{algorithmic}

\usepackage{url}
\usepackage{graphicx, tikz} % Useful for graphics
\usetikzlibrary{calc,positioning,shapes,fit,arrows,shadows}
\tikzstyle{line} = [ draw, -latex']

\usepackage{amsthm}

\newtheorem{prop}{Proposition}[section]

\newtheorem{obs}{Observation}[section]

\allowdisplaybreaks

\usepackage{rotating}
\usepackage{pdflscape}

\graphicspath{{figures/}} % Graphics will be here

\title{New Formulation and Strong MISOCP Relaxations \\ for AC Optimal Transmission Switching Problem}% \footnote{Under second review for {\it IEEE Transactions on Power Systems}, invited for a special issue}} 

\author{Burak Kocuk, Santanu S. Dey, X. Andy Sun 
}

\DeclareMathOperator{\atan2}{atan2}

\begin{document}
\maketitle

\begin{abstract}
As the modern transmission control and relay technologies evolve, transmission line switching has become an important option in power system operators' toolkits {to reduce operational cost and improve system reliability.} %since it promises cost reduction and better reliability.  
Most recent research has {relied on the DC approximation of the power flow model in the optimal transmission switching problem}. However, it is  known that DC approximation may {lead to} inaccurate flow solutions and also overlook stability {issues}. In this paper, we focus on the {optimal transmission switching problem with the full AC power flow model}, abbreviated as AC OTS. We propose  {a new exact} formulation for {AC OTS} and its {mixed-integer second-order conic programming (MISOCP)} relaxation. We improve this relaxation via several types of {strong} valid inequalities inspired by the recent development for the closely related AC Optimal Power Flow (AC OPF) problem \cite{kocuk2015}. We also propose a practical algorithm to obtain high quality feasible solutions {for the AC OTS problem}. {Extensive computational experiments show that the proposed formulation and algorithms efficiently solve IEEE standard and congested instances and lead to significant cost benefits with provably tight bounds}. 
\end{abstract}

\section{Introduction}

{ Transmission switching, as an emerging operational scheme, has gained considerable attention in both industry and academia in the recent years \cite{oneill.et.al:05, Fisher,  Hedman09,  Khodaei10, Hedman11}. Switching on and off transmission lines, therefore, changing the network topology in the real-time operation, may bring several benefits that the traditional economic dispatch cannot offer, such as 
reducing the total {operational} cost{\cite{Fisher, Hedman2010, Han2016}}, 
mitigating transmission congestion{\cite{Thompson2009}},  
clearing contingencies{\cite{Khanabadi2013, Korad2013}},  and
{improving do-not-exceed limits\cite{Korad2015}}
.}

Previous literature on OTS mainly {relies on the DC approximation of the power flow model to avoid} the mathematical complexity induced by the non-convexity of AC power flow equations {(see e.g. \cite{oneill.et.al:05, Fisher, Ruiz2011, Ruiz2012})}. This {DC} version of the {OTS problem} can be modeled as a  mixed-integer linear program (MILP), {which is a} computationally {challenging problem} and several heuristic method are proposed % to explore high quality feasible solutions 
\cite{Barrows12,Fuller12,Wu13}. In a recent work \cite{kocuk2014Switch}, {the authors} propose a {new formulation and a} class of valid inequalities to {exactly} solve the MILP problem. However, even if this problem can be solved quickly, it has been recognized that the optimal topology obtained by solving DC transmission switching is not guaranteed to be AC feasible, also it may over-estimate cost improvements {and overlook stability issues} \cite{hijazi2013}. 

The AC optimal transmission switching problem (AC OTS) is {much less explored}. In \cite{hijazi2013}, a convex relaxation of AC OTS is proposed based on trigonometric outer-approximation. The problem is formulated as a mixed integer nonlinear program (MINLP) and solved using the solver BONMIN to obtain upper bounds. In \cite{Soroush14}, a new ranking heuristic is proposed based on the economic dispatch solutions and the corresponding dual variables. In \cite{barrow.blumsack.hines:14},  {DC OTS} is solved {first and then} a heuristic correction mechanism is utilized to restore AC feasibility of the solutions. In this paper, we aim to push the control scheme for transmission switching closer to the real-world power system operation by proposing a new {exact} formulation and an efficient algorithm for the AC OTS problem. % in a more rigorous level. 

{There are several closely related problems in the literature, which  involve line switching decisions, such as the network configuration problem \cite{Jabr12, Ferreira14}, transmission system planning \cite{Jabr13}, and intentional islanding \cite{Trodden14}. The main ideas of these works are based on conic relaxations or piecewise linear approximations of the non-convex power flow equations. }

Our study starts from the recent advances in a related fundamental problem in power system analysis, namely the AC Optimal Power Flow (AC OPF) problem, which minimizes the generation cost to satisfy load and various physical {constraints} represented in the AC power flow constraints, {while} the power network topology is kept unchanged. It is {demonstrated by several authors} that convex relaxations, especially semidefinite programming (SDP) relaxations, of the AC OPF problem {provide tight lower bounds} on {standard} IEEE test instances \cite{bai08, lavaei12, molzahn2013, molzahn2014}. However, the computational burden of solving large-scale SDP relaxations is still unwieldy. {To solve for large-scale systems}, one {may need to} turn to computationally less demanding alternatives such as quadratic convex \cite{coffrin2014, hijazi2013, coffrin2015} or linear programming relaxations \cite{bienstock2014}.

In a recent work \cite{kocuk2015}, we proposed several strong second-order cone programming (SOCP) relaxations for AC OPF, which produce extremely high quality feasible AC solutions (not dominated by the SDP relaxations) in a time that is an order of magnitude faster than solving the SDP relaxations. {In this paper, we extend these new techniques to the more challenging AC OTS problem}. In particular, we first formulate the AC OTS problem as an {MINLP} problem. Then, we propose a mixed-integer second-order cone programming (MISOCP) relaxation, which relaxes the non-convex AC power flow constraints to a set of convex quadratic constraints, represented in the form of SOCP constraints. The paper then provides several techniques to strengthen this MISOCP relaxation by adding several types of valid inequalities. Some of these valid inequalities have demonstrated {to have} excellent performance for the AC OPF in \cite{kocuk2015}, and some others are specifically developed for the AC OTS problem. Finally, we also propose practical algorithms that utilize the solutions from the MISOCP relaxation to obtain high quality feasible solutions for the AC OTS problem.

The rest of the paper is organized as follows: In Section \ref{sec:opf} we formally define {AC} OPF and present two {exact} formulations. In Section \ref{sec:ots}, we present {AC OTS} as an MINLP problem and discuss its MISOCP relaxation. Then, we propose several valid inequalities in Section \ref{sec:valid ineq} and develop a practical algorithm to solve AC OTS in Section \ref{sec:alg}. We present the results of our extensive computational experiments in Section \ref{sec:comp expr}. Finally, some concluding remarks are given in Section \ref{sec:concl}.

\section{AC Optimal Power Flow}
\label{sec:opf}

Consider a power network $\mathcal{N} = (\mathcal{B},\mathcal{L})$, where $\mathcal{B}$ and $\mathcal{L}$  respectively denote the set of buses and transmission lines. Generation units are connected to a subset of buses, denoted as $\mathcal{G}\subseteq \mathcal{B}$. 
%We assume that there is load at every bus. 
The aim of the AC optimal power flow (OPF) problem is to satisfy demand at all buses with the minimum total production costs of generators such that the solution obeys the physical laws (e.g., Ohm's  and Kirchoff's Law) and other operational restrictions (e.g., transmission line flow limit constraints).

Let $Y \in \mathbb{C}^{|\mathcal{B}| \times |\mathcal{B}|}$ denote the nodal admittance matrix, which has components $Y_{ij}=G_{ij} + \mathrm{i}B_{ij}$ for each line $(i,j)\in\mathcal{L}$, and $G_{ii}=g_{ii}-\sum_{j\ne i} G_{ij}, B_{ii}=b_{ii}-\sum_{j\ne i} B_{ij}$, where $g_{ii}$ (resp. $b_{ii}$) is the shunt conductance (resp. susceptance) at bus $i\in\mathcal{B}$ and $\mathrm{i}=\sqrt{-1}$. 
Let $p_i^g, q_i^g$ (resp. $p_i^d, q_i^d$) be the real and reactive power output of the generator (resp. load) at bus $i$. The complex voltage $V_i$ at bus $i$ can be expressed either in the rectangular form as $V_i = e_i+\mathrm{i} f_i$ or in the polar form as $V_i = |V_i|(\cos\theta_i+\mathrm{i}\sin\theta_i)$, where $|V_i|=\sqrt{e_i^2 + f_i^2}$ is the voltage magnitude and $\theta_i$ is the phase angle. Real and reactive power on line $(i,j)$ are denoted by  $p_{ij}$ and $q_{ij}$, respectively and computed as follows:
\begin{equation}\label{flow definitions}
\begin{split}
p_{ij} &= -G_{ij}(e_i^2+f_i^2) + G_{ij}(e_ie_j+f_if_j) -B_{ij}(e_if_j-e_jf_i) \\
q_{ij} &= \ \ \,B_{ij}(e_i^2+f_i^2) -  B_{ij}(e_ie_j+f_if_j) - G_{ij}(e_if_j-e_jf_i).
\end{split}
\end{equation}
%where $\delta(i)$ denotes the buses directly connected to $i$ by a transmission line.

With the above notation, the AC OPF problem is given in the so-called rectangular formulation as follows:
\begin{subequations} \label{rect}
\begin{align}
  \min  &\hspace{0.5em}  \sum_{i \in \mathcal{G}} C_i(p_i^g)  \label{obj} \\
  \mathrm{s.t.}   &\hspace{0.5em} p_i^g-p_i^d = g_{ii}(e_i^2+f_i^2) + \sum_{j \in \delta(i)}p_{ij}   & i& \in \mathcal{B} \label{activeAtBus} \\
  & \hspace{0.5em} q_i^g-q_i^d = -b_{ii}(e_i^2+f_i^2) + \sum_{j \in \delta(i)}q_{ij}  & i& \in \mathcal{B} \label{reactiveAtBus} \\
  & \hspace{0.5em} \underline V_i^2 \le e_i^2+f_i^2 \le \overline V_i^2    & i& \in \mathcal{B} \label{voltageAtBus} \\
  & \hspace{0.5em}  p_{ij}^2 + q_{ij}^2   \le (S_{ij}^{\text{max}})^2  &(&i,j) \in \mathcal{L} \label{powerOnArc} \\
  & \hspace{0.5em}  p_i^{\text{min}}  \le p_i^g \le p_i^{\text{max}}     & i& \in \mathcal{G} \label{activeAtGenerator} \\
  & \hspace{0.5em} q_i^{\text{min}}  \le q_i^g \le q_i^{\text{max}}     & i& \in \mathcal{G} \label{reactiveAtGenerator}, \\
\notag  & \hspace{0.5em} \eqref{flow definitions}.
\end{align}
\end{subequations}
The objective function $C_i(p_i^g)$ is typically linear or convex quadratic in the real power output $p_i^g$ of generator $i$. Constraints \eqref{activeAtBus} and \eqref{reactiveAtBus} correspond to the conservation of active and reactive power flows at each bus, respectively. Here, $\delta(i)$ denotes the set of neighbor buses of bus $i$. Constraint \eqref{voltageAtBus} restricts voltage magnitude at each bus.
%As noted above, $\underline V_i$ and $\overline V_i$ are both close to  1 p.u. at each bus $i$. 
 Constraint \eqref{powerOnArc} puts an upper bound on the apparent power on each line.
Finally, constraints \eqref{activeAtGenerator} and \eqref{reactiveAtGenerator}, respectively, limit the active and reactive power output of each generator to respect its physical capability.

Note that the rectangular formulation \eqref{rect} is a non-convex quadratic optimization problem. However, we note that all the nonlinearity and non-convexity comes from one of the following three forms:
(1) $e_i^2+f_i^2=|V_i|^2$, (2) $e_ie_j+f_if_j=|V_i||V_j|\cos(\theta_i-\theta_j)$, (3) $e_if_j-f_ie_j=-|V_i||V_j|\sin(\theta_i-\theta_j)$. We define new variables $c_{ii}$, $c_{ij}$ and $s_{ij}$ for each bus $i$ and each transmission line $(i,j)$ to capture the {non-convexity}. In particular, we {define for each $i\in\mathcal{B}$ and $(i,j)\in\mathcal{L}$},
\begin{equation}\label{cs definitions}
c_{ii} := e_i^2 + f_i^2, \quad c_{ij} := e_ie_j+f_if_j, \quad s_{ij} :=e_if_j-e_jf_i. 
\end{equation}
Now, we introduce an {equivalent,} alternative formulation of the OPF problem as follows:
\begin{subequations} \label{SOCP}
\begin{align}
 \min  &\hspace{0.5em}  \sum_{i \in \mathcal{G}} C_i(p_i^g) \\
  \mathrm{s.t.}   &\hspace{0.5em} p_i^g-p_i^d = g_{ii}c_{ii} + \sum_{j \in \delta(i)} p_{ij}  & i& \in \mathcal{B} \label{activeAtBusR} \\
  & \hspace{0.5em} q_i^g-q_i^d = -b_{ii}c_{ii} + \sum_{j \in \delta(i)} q_{ij} & i& \in \mathcal{B} \label{reactiveAtBusR} \\ 
  & \hspace{0.5em} p_{ij} = -G_{ij}c_{ii}+ G_{ij}c_{ij}-B_{ij}s_{ij}  &(&i,j) \in \mathcal{L} \\
  & \hspace{0.5em} q_{ij} = \ \ \,B_{ij}c_{ii} -  B_{ij}c_{ij} - G_{ij}s_{ij} &(&i,j) \in \mathcal{L} \\
  & \hspace{0.5em} \underline V_i^2 \le c_{ii} \le \overline V_i^2    & i& \in \mathcal{B}  \label{voltageAtBusR} \\
  & \hspace{0.5em} c_{ij}=c_{ji}, \ \ s_{ij}=-s_{ji}    &(&i,j) \in \mathcal{L} \label{cosine_sine}\\
  & \hspace{0.5em}  c_{ij}^2+s_{ij}^2  = c_{ii}c_{jj}     &(&i,j) \in \mathcal{L} \label{coupling}\\
  & \hspace{0.5em}  \theta_j - \theta_i = \atan2({s_{ij}},{c_{ij}})  &(&i,j) \in \mathcal{L},  \label{arctan cons} \\
\notag  & \hspace{0.5em}  \eqref{powerOnArc}\text{-}\eqref{reactiveAtGenerator}. 
\end{align}
\end{subequations}
A variant of this formulation without \eqref{arctan cons} was previously proposed in  \cite{Exposito99} and \cite{Jabr06} for radial networks {(also see \cite{kocuk2014})}
while it was later adapted to general networks in \cite{Jabr07,Jabr08}. 

\section{AC Optimal Transmission Switching}
\label{sec:ots}

{AC} Optimal Transmission Switching (AC OTS) is a variant of {the AC OPF problem} in which transmission lines are allowed to be switched on and off to reduce the total cost of dispatch. %This counterintuitive phenomenon is due to a form of Braess' paradox. 
AC OTS can be formulated as an optimization problem, which aims to find a topology with the least cost while {achieving feasible AC power flow solutions}. In this section, we first formulate AC OTS as an MINLP and then, propose an MISOCP relaxation to obtain lower bounds. {We will use OTS (resp. OPF) to denote AC OTS (resp. AC OPF) for brevity, unless stated otherwise}.

\subsection{MINLP Formulation}
Mathematical programming formulation of OTS can be {stated} with the same variables as used in OPF with the addition of a set of binary variables, denoted by $x_{ij}$, for each line. The variable {$x_{ij}$} takes the value one if the corresponding line {$(i,j)$} is switched on, and zero otherwise. Then, OTS is formulated as {the following MINLP problem}:
\begin{subequations}\label{OTS reform}
\begin{align}
 \min &\hspace{0.0em}  \sum_{i \in \mathcal{G}} C_i(p_i^g) \\
  \mathrm{s.t. }  \   & \hspace{0.0em}p_{ij} = (-G_{ij}c_{ii} + G_{ij}c_{ij} -B_{ij}s_{ij}) x_{ij}  &(&i,j) \in \mathcal{L}  \label{realOnArc x} \\
  & \hspace{0.0em}q_{ij} = \ \ \,(B_{ij}c_{ii} -  B_{ij}c_{ij} - G_{ij}s_{ij})x_{ij}&(&i,j) \in \mathcal{L}   \label{reactiveOnArc x} \\
  & \hspace{0.0em}  (c_{ij}^2+s_{ij}^2  - c_{ii}c_{jj} )x_{ij} = 0     &(&i,j) \in \mathcal{L} \label{coupling binary}\\
  & \hspace{0.0em} (\theta_j- \theta_i - \atan2(s_{ij},c_{ij}))x_{ij} = 0     &(&i,j) \in \mathcal{L}   \label{arctan cons binary}  \\
 & \hspace{0.0em} x_{ij}  \in\{ 0,1\} &(&i,j) \in \mathcal{L},  \label{binary x}\\
\notag  & \hspace{0.0em}  \eqref{powerOnArc}\text{-}\eqref{reactiveAtGenerator}, \eqref{activeAtBusR}\text{-}  \eqref{reactiveAtBusR}, \eqref{voltageAtBusR}\text{-} \eqref{cosine_sine}.
\end{align}
\end{subequations}
Here, constraints \eqref{realOnArc x} and \eqref{reactiveOnArc x} guarantee that real and reactive flow on every line takes the associated values if the line is
switched on and zero otherwise. Similarly,  constraints \eqref{coupling binary} and  \eqref{arctan cons binary}    are active only when the corresponding binary variable takes the value one.

{We also note that  the model \eqref{OTS reform} can be appropriately modified to include circuit breakers between bus bars \cite{VanAcker}.}

\subsection{MISOCP Relaxation}
Now, we propose an MISOCP relaxation of OTS \eqref{OTS reform}. For notational convenience, let $ \underline c_{ii} = \underline V_i^2$ and $ \overline c_{ii} =\overline V_i^2$. Here, we extend the definition of variables $c_{ij}$ and $s_{ij}$, which now take the values as before when the corresponding line is switched on and zero otherwise. We also denote lower and upper bounds {of} $c_{ij}$ (resp. $s_{ij}$) {as} $\underline c_{ij}$ (resp. $\underline s_{ij}$) and $\overline c_{ij}$ (resp. $\overline s_{ij}$), respectively, when the line is switched on. Next, we define new variables $c_{ii}^j:=c_{ii}x_{ij}$. Using this notation, we present an MISOCP relaxation as follows:
\begin{subequations} \label{misocp}
\begin{align}
  \min  &\hspace{0.5em}  \sum_{i \in \mathcal{G}} C_i(p_i^g)   \\
  \mathrm{s.t.}   
  & \hspace{0.5em}p_{ij} = -G_{ij}c_{ii}^j + G_{ij}c_{ij} -B_{ij}s_{ij}     &(&i,j) \in \mathcal{L} \label{realOnArc 2}\\
  & \hspace{0.5em}q_{ij} = \ \ \,B_{ij}c_{ii}^j -  B_{ij}c_{ij} - G_{ij}s_{ij}&(&i,j) \in \mathcal{L}\label{reactiveOnArc 2} \\
  & \hspace{0.5em} \underline c_{ij}x_{ij}\le c_{ij} \le \overline c_{ij} x_{ij}   &(&i,j) \in \mathcal{L} \label{c dummy} \\
  & \hspace{0.5em} \underline s_{ij}x_{ij}\le s_{ij} \le \overline s_{ij} x_{ij}   &(&i,j) \in \mathcal{L} \label{s dummy} \\
  & \hspace{0.5em}  \underline c_{ii}x_{ij} \le c_{ii}^j \le \overline c_{ii} x_{ij}   &(&i,j) \in \mathcal{L} \label{c_dummy} \\
  & \hspace{0.5em}  c_{ii}- \overline c_{ii}(1-x_{ij}) \le c_{ii}^j  &(&i,j) \in \mathcal{L} \label{c_dummy 1}  \\
  & \hspace{0.5em}  c_{ii}^j \le c_{ii} - \underline  c_{ii}(1- x_{ij})   &(&i,j) \in \mathcal{L} \label{c_dummy 2}  \\
  & \hspace{0.5em}   c_{ij}^2+s_{ij}^2  \le c_{ii}^jc_{jj}^i     &(&i,j) \in \mathcal{L}, \label{coupling new} \\
\notag  & \hspace{0.5em}  \eqref{powerOnArc}\text{-}\eqref{reactiveAtGenerator}, \eqref{activeAtBusR}\text{-} \eqref{reactiveAtBusR}, \eqref{voltageAtBusR}\text{-} \eqref{cosine_sine}, \eqref{binary x}.
\end{align}
\end{subequations}
Here, constraints \eqref{realOnArc 2} and \eqref{reactiveOnArc 2} again guarantee that flow variables takes the correct value when the line is switched on and zero otherwise, due to constraints \eqref{c dummy}-\eqref{c_dummy}. On the other hand,  \eqref{c_dummy 1}-\eqref{c_dummy 2} restrict that $c_{ii}^j$ takes value $c_{ii}$ when line in switched on. We note that constraints \eqref{c_dummy}-\eqref{c_dummy 2} are precisely the  McCormicks
envelopes \cite{mccormick} applied to $c_{ii}^j=c_{ii}x_{ij}$. Finally,  \eqref{coupling new} is the SOCP relaxation of \eqref{coupling binary}.

We note that the non-convex constraint \eqref{arctan cons binary} is dropped altogether to obtain the MISOCP relaxation \eqref{misocp}. In the next section, we propose three ways to incorporate the constraint \eqref{arctan cons binary} back into the MISOCP relaxation.

\section{Valid Inequalities}
\label{sec:valid ineq} 

In this section, we propose three methods to strengthen the MISOCP relaxation \eqref{misocp}. They are based on the strengthening methods we recently proposed for the SOCP relaxation of the AC OPF problem in \cite{kocuk2015}, which are combined with integer programming techniques.
In Section \ref{sec:arctan envelope x}, we construct a polyhedral envelope for the arctangent constraint  \eqref{arctan cons binary} in 3-dimension.
In Section \ref{SOCP SDP sep x}, we propose a disjunctive cut generation scheme    that separates a given SOCP solution from the SDP cones. 
%This approach takes the advantage of the efficiency of the SOCP relaxation and the accuracy of the SDP relaxation.
In Section \ref{McCormick sep x}, we propose another disjunctive cut generation scheme  that separates a given SOCP solution from a newly-proposed cycle based McCormick relaxation of the OPF problem. 
Finally, in Section \ref{sec:bounding x}, we propose variable bounding techniques that provide tight variable bounds, which is essential for the success of the proposed approach.

\subsection{Arctangent Envelopes}
\label{sec:arctan envelope x}

First, we propose a convex outer-approximation of the angle condition \eqref{arctan cons binary} to the MISOCP relaxation.
%An alternative approach to incorporate cycles would be to use  the arctangent constraint \eqref{arctan cons} directly. In this section, we propose convex envelopes for this particular constraint.
Our construction uses four linear inequalities to approximate the convex envelope for the following set defined by the arctangent constraint \eqref{arctan cons binary} for each line $(i,j)\in\mathcal{L}$, %When the lower bound $\underline c_{ij}$ on $c_{ij}$ is positive, \eqref{arctan cons} is equivalent to the following constraint,}
\begin{equation}
\mathcal{AT}:=\left\{(c,s,\theta)\in\mathbb{R}^3 : \theta = \arctan \left ( \frac sc \right), (c,s) \in B \right\},	\label{arctan constraint}
\end{equation}
where we denote $\theta= \theta_j - \theta_i$ and drop $(i,j)$ indices for brevity and define the box $B := [\underline c, \overline c]\times[\underline s, \overline s]$. We also assume $\underline c>0$. The four corners of the box correspond to four points in the $(c,s,\theta)$ space: 
%To start with, let us denote  the corners of the box $[\underline c, \overline c]\times[\underline s, \overline s]$ as 
\begin{equation}
\begin{split}
   z^1 &= (\underline c, \overline s, \arctan \left( {\overline s}/{\underline c} \right) ), \quad z^2 = (\overline c, \overline s, \arctan \left( {\overline s}/{\overline c} \right) ),\\ 
   z^3 &= (\overline c,  \underline s, \arctan \left( {\underline s}/{\overline c} \right)), \quad z^4 = (\underline c, \underline s, \arctan \left( {\underline s}/{\underline c} \right)). 
\end{split}
\end{equation}

Let us first focus on the upper envelopes. {Proposition \ref{prop:arctan x lower} is adapted from \cite{kocuk2015} to the case of OTS:}
\begin{prop}\label{prop:arctan x lower}
Let $\theta = \gamma_1 + \alpha_1 c + \beta_1 s$  and $\theta = \gamma_2 + \alpha_2 c + \beta_2 s$ be the  planes passing through points $\{z^1,z^2, z^3\}$,  and $\{z^1,z^3, z^4\}$, respectively.
Then, for $k=1,2$, we have
\begin{equation}\label{eq:arctan-upper}
\gamma_k' + \alpha_k c + \beta_k s + (2\pi-\gamma_k' )(1-x) \ge \arctan \left( \frac{s}{c} \right)\end{equation}  for all $(c,s) \in B$
 with $\gamma_k' = \gamma_k + \Delta \gamma_k$ where
\begin{equation} \label{error problem}
\Delta \gamma_k =   \max_{(c,s) \in B} \left \{  \arctan \left( \frac{s}{c} \right) - (\gamma_k+ \alpha_k c + \beta_k s )  \right \}. 
\end{equation} 
\end{prop}
%Note that by the construction of \eqref{error problem}, it is evident that $\gamma_k'+\alpha_k c + \beta_k s$ dominates  $\arctan(s/c)$ over the box $B$ when $x=1$ and the inequalities become redundant otherwise.
The nonconvex optimization problem  \eqref{error problem} can be solved by enumerating all possible Karush-Kuhn-Tucker (KKT) points.
A similar argument can be used to construct lower envelopes as well. See \cite{kocuk2015} for details.

\subsection{SDP Disjunction}
\label{SOCP SDP sep x}

In the second method to strengthen the MISOCP relaxation \eqref{misocp}, we propose a cutting plane approach to separate a given SOCP relaxation solution from the feasible region of the SDP relaxation of cycles. To start with, let us consider a cycle with the set of lines $C$ and the set of buses  $\mathcal{B}_C$.
Let $v\in \mathbb{R}^{2|C|}$ be a vector of bus voltages defined as $v=[e; f]$ such that $v_i=e_i$ for $i\in\mathcal{B}$ and $v_{i'}=f_i$ for $i'=i+|C|$. Observe that if we have a set of $c,s$ variables satisfying the definitions in \eqref{cs definitions} and a matrix variable $W=vv^T$, then the following  relationship holds between $c,s,x$ and $W$,
\begin{subequations}\label{eq:sdplin-x}
\begin{align}
&c_{ij} = (W_{ij} + W_{i'j'})x_{ij} & (&i,j)\in C \\ 
&s_{ij} =(W_{ij'} - W_{ji'})x_{ij} & (&i,j)\in C \\ 
&c_{ii} =  W_{ii} + W_{i'i'}   & i&\in\mathcal{B}_C \label{dummy1} \\
&\underline c_{ij}x_{ij}\le c_{ij} \le \overline c_{ij}x_{ij} & (&i,j)\in C \label{dummy2}   \\ 
&\underline s_{ij}x_{ij} \le s_{ij} \le \overline s_{ij}x_{ij} & (&i,j)\in C \label{dummy3}  \\
&\underline c_{ii}\le c_{ii} \le \overline c_{ii}  & i&\in\mathcal{B}_C \label{dummy4} \\
&c_{ii}^j = c_{ii}x_{ij} & (&i,j)\in C \\ 
&x_{ij} \in\{0,1\} & (&i,j)\in C \label{dummy5} \\
& W \succeq 0.\label{dummy6} 
\end{align}
\end{subequations}
Let us define $\mathcal{S}:=\{(c,s,x) : \exists W: \eqref{eq:sdplin-x}\}$. Clearly, any feasible solution to the OTS formulation \eqref{OTS reform} must also satisfy \eqref{eq:sdplin-x}. Therefore, any valid inequality for $\mathcal{S}$ is also valid for the formulation  \eqref{OTS reform}.

Note that $\mathcal{S}$ is a mixed-integer set. Ideally, one would be interested in finding conv$(\mathcal{S})$ to generate strong valid inequalities. {However, this is a quite computationally challenging task, no easier than solving the original MINLP.} Instead, we outer-approximate  conv$(\mathcal{S})$ and obtain cutting planes by utilizing a simple disjunction for a cycle $C$: Either every line is active, that is $ \sum_{(i,j)\in C} x_{ij} = |C|$, or at least one line is disconnected, that is $\sum_{(i,j)\in C} x_{ij} \le |C| - 1$. Below, we approximate these two disjunctions.

\underline{Disjunction 1}: In the first disjunction, we have $x_{ij}=1$ for all $(i,j) \in C$. Let us consider the following constraints
\begin{subequations}\label{sdp disj1}
\begin{align}
&c_{ij}  = W_{ij} + W_{i'j'} & (i,j)\in C\\ 
&s_{ij} =  W_{ij'} - W_{ji'} & (i,j)\in C\\ 
  & c_{ii}= c_{ii}^j   &(i,j) \in C \label{sdp disj1_} \\
&x_{ij} = 1 & (i,j)\in C \label{binary1},
\end{align}
\end{subequations}
and define $\mathcal{S}_1 := \{(c,s,x) : \exists W : \eqref{sdp disj1}, \eqref{dummy1}-\eqref{dummy4}, \eqref{dummy6} \}$.

\underline{Disjunction 0}: In the second disjunction,  $x_{ij}=0$ for some $(i,j) \in C$. Let us consider the following constraints
\begin{subequations}\label{sdp disj0}
\begin{align}
& c_{ij}^2+s_{ij}^2  \le c_{ii}^jc_{jj}^i   & (i,j)\in C\\ 
& \underline c_{ii}x_{ij} \le c_{ii}^j \le \overline c_{ii} x_{ij}   &(i,j) \in C \\
  & c_{ii}- \overline c_{ii}(1-x_{ij}) \le c_{ii}^j   &(i,j) \in C \\
  & c_{ii}^j \le c_{ii} - \underline  c_{ii}(1- x_{ij})   &(i,j) \in C \\
&0 \le x_{ij} \le 1 & (i,j)\in C\\ 
& \sum_{(i,j) \in C} x_{ij} \le |C| - 1,
\end{align}
\end{subequations}
and define $\mathcal{S}_0 := \{(c,s,x)  : \eqref{sdp disj0}, \eqref{dummy2}\text{-}\eqref{dummy4} \}$.

We note that  both $\mathcal{S}_1 $  and $\mathcal{S}_0$ are conic representable. In particular, these bounded sets are respectively semidefinite and second-order cone representable. Therefore, conv$(\mathcal{S}_1 \cup \mathcal{S}_0)$ is also conic representable (see  Appendix \ref{app:union}  on how to obtain a representation as an extended formulation), and by construction, contains $\mathcal{S}$.  

Now, suppose a point $(c^*,s^*,x^*)$ is given. We want to decide whether this point belongs to conv$(\mathcal{S}_1 \cup \mathcal{S}_0)$ or otherwise, find a separating hyperplane. Given that we have an extended semidefinite representation for conv$(\mathcal{S}_1 \cup \mathcal{S}_0)$, we can solve an SDP separation problem to achieve this. See Appendix \ref{app:sep}.

\subsection{McCormick Disjunction}
\label{McCormick sep x}

The last method to strengthen the MISOCP relaxation \eqref{misocp} is based on a new cycle-based OPF formulation we propose in \cite{kocuk2015}. 
The key observation is as follows: instead of satisfying the angle condition \eqref{arctan cons binary} for each $(i,j) \in \mathcal{L}$, we guarantee that angle differences sum up to 0 modulo $2\pi$ over every cycle $C$ in the power network if all the lines of the cycle $C$ are switched on, i.e.
\begin{equation}
\bigl(\sum_{(i,j)\in C } \theta_{ij} - 2\pi k\bigr)\prod_{(i,j)\in C }x_{ij} = 0, \quad \text{ for some } k \in \mathbb{Z}, \label{angle = 0.} 
\end{equation}
where $\theta_{ij} := \theta_j-\theta_i$. 

Next, we consider
\begin{subequations}\label{eq:mcclin-x}
\begin{align}
& \bigl [ \cos \bigl ( \sum_{(i,j)\in C } \theta_{ij}  \bigr ) - 1 \bigr] \prod_{(i,j)\in C }x_{ij} = 0 \label{eq:mcclin-x1} \\
& c_{ij} = \sqrt{c_{ii}c_{jj}}\cos \theta_{ij} x_{ij} \qquad \qquad\qquad\qquad (i,j) \in C\label{eq:mcclin-x2} \\
& s_{ij} = \sqrt{c_{ii}c_{jj}}\sin \theta_{ij} x_{ij} \qquad \qquad\qquad\qquad (i,j) \in C, \label{eq:mcclin-x3}\\
&\notag \eqref{dummy2}-\eqref{dummy5}.
%&\underline c_{ij}x_{ij}\le c_{ij} \le \overline c_{ij}x_{ij} & (i,j)\in C \label{dummy2}   \\ 
%&\underline s_{ij}x_{ij} \le s_{ij} \le \overline s_{ij}x_{ij} & (i,j)\in C \label{dummy3}  \\
%&\underline c_{ii}\le c_{ii} \le \overline c_{ii}  & i\in\mathcal{B}_C \label{dummy4} \\
%&c_{ii}^j = c_{ii}x_{ij} & (i,j)\in C \\ 
%&x_{ij} \in\{0,1\} & (i,j)\in C.
\end{align}
\end{subequations}
%where $i'=i+|{\mathcal B}|$ and $j'=j+|{\mathcal B}|$. %Here we suppress labeling any variable with $C$ for brevity. It is understood that such construction is done for each cycle in a cycle basis. 
Here, \eqref{eq:mcclin-x1} is equivalent to \eqref{angle = 0.} and \eqref{eq:mcclin-x2}-\eqref{eq:mcclin-x3} follow from the definition of $c,s$ variables.
Let us define $\mathcal{M}:=\{(c,s,x) : \exists \theta: \eqref{eq:mcclin-x}, {\eqref{dummy2}-\eqref{dummy5}}\}$. Again, observe that any feasible solution to the OTS formulation \eqref{OTS reform} must also satisfy \eqref{eq:mcclin-x}. Therefore, any valid inequality for $\mathcal{M}$ is also valid for the formulation  \eqref{OTS reform}.

%Note that $\mathcal{S}$ is a mixed-integer set. Ideally, one would be interested in finding conv$(\mathcal{S})$ to generate strong valid inequalities. Instead, we ``outer-approximate"  conv$(\mathcal{S})$ and obtain cutting planes by utilizing a simple disjunction for a cycle $C$: Either every line is active, that is $ \sum_{(i,j)\in C} x_{ij} = |C|$, or at least one line is disconnected, that is $\sum_{(i,j)\in C} x_{ij} \le |C| - 1$. Below, we approximate these two disjunctions:

We again follow a similar procedure to the previous section and consider  two disjunctions for a  cycle $C$.

\underline{Disjunction 1}: In the first disjunction, we have $x_{ij}=1$ for all $(i,j) \in C$. Note that  \eqref{eq:mcclin-x1}  reduces to 
\begin{equation*}
\cos \bigl ( \sum_{(i,j)\in C } \theta_{ij}  \bigr ) = 1.
\end{equation*}
Now, we can expand the cosine appropriately and replace  $\cos(\theta_{ij})$'s and $\sin(\theta_{ij})$'s in terms of $c,s$ variables following \eqref{eq:mcclin-x2}-\eqref{eq:mcclin-x3}. This transformation yields a homogeneous polynomial, denoted by $p_C$, in terms of only $c,s$ variables, and an equivalent constraint $p_C = 0$. However, $p_C$ can have up to $2^{|C|-1}+1$ monomials and each monomial of degree $|C|$. In  \cite{kocuk2015}, we propose  a  method, which is used to ``bilinearize" this high degree polynomial by decomposing larger cycles into smaller ones by the addition of artificial lines and corresponding variables. We refer the reader to \cite{kocuk2015} for details.

Using the proposed decomposition scheme, we obtain a set of bilinear polynomials
$
q_k(c,s, \tilde c, \tilde s) = 0, \ k \in \mathcal{K}_C,
$
for a given cycle $C$. Here, $\tilde c, \tilde s$ denote the extra variables used in the construction. %, which represent the same quantities for artificial lines.

Finally, we use McCormick envelopes for each bilinear constraint to linearize the system of polynomials. For a given cycle $C$, consider
 the McCormick relaxation of the bilinear cycle constraints, which   can be written compactly as follows:
\vspace{-0mm}
\begin{subequations}\label{mcc disj1}
\begin{align}
& Az + \tilde A \tilde z + B y \le c \label{mccormick and bounds} \\
& Ey = 0. \label{cycle equalities} 
\end{align}
\end{subequations}
Here, $z$ is a vector composed of the $c,s$ variables, $\tilde z$ is a vector composed of the additional $\tilde c, \tilde s$ variables introduced in the cycle decomposition, and $y$ is a vector of new variables defined to linearize the bilinear terms in the cycle constraints. Constraint \eqref{mccormick and bounds} contains the McCormick envelopes of the bilinear terms and bounds  on the $c,s$ variables, while \eqref{cycle equalities} includes the linearized cycle equality constraints. Finally, we define the set $\mathcal{M}_1 := \{(c,s,x) : \exists (\tilde c, \tilde s) : \eqref{mcc disj1}, \eqref{dummy2}\text{-}\eqref{dummy4}, \eqref{sdp disj1_}\text{-}\eqref{binary1} \}$.

\underline{Disjunction 0}: In  the second disjunction,  $x_{ij}=0$ for some $(i,j) \in C$. We take $\mathcal{M}_0 := \mathcal{S}_0$.

We note that  both $\mathcal{M}_1 $  and $\mathcal{M}_0$ are conic representable. In particular, these bounded sets are respectively polyhedral and second-order cone representable. Therefore, conv$(\mathcal{M}_1 \cup \mathcal{M}_0)$ is also conic representable, and by construction, contains $\mathcal{M}$.  

Now, suppose a point $(c^*,s^*,x^*)$ is given. We want to decide whether this point belongs to conv$(\mathcal{M}_1 \cup \mathcal{M}_0)$ or otherwise, find a separating hyperplane. Given that we have an extended second-order cone representation for conv$(\mathcal{M}_1 \cup \mathcal{M}_0)$, we can solve an SOCP separation problem.

In our computations, we observed that stronger cuts are obtained by combining SDP and McCormick Disjunction. In particular, we separate cutting planes from
  conv$( (\mathcal{S}_1 \cap \mathcal{M}_1) \cup \mathcal{S}_0)$ by solving SDP separation problems.

\subsection{Obtaining Variable Bounds}
\label{sec:bounding x}

Note that the  arctangent envelopes and the McCormick relaxations are more effective  when tight variable  upper/lower bounds are available for the $c$ and $s$ variables. Now, we explain how we obtain good bounds {for these variables}, which is the key ingredient in the success of our proposed methods.
 
Observe that $c_{ij}$ and $s_{ij}$  do not have explicit variable bounds except the implied bounds due to (\ref{voltageAtBusR})  and (\ref{coupling}) as
\begin{equation*}
-\overline V_i \overline V_j \le c_{ij}, s_{ij} \le \overline V_i \overline V_j \quad (i,j) \in \mathcal{L}.
\end{equation*}
However, these bounds {may be quite}  loose, {especially when the} phase angle differences are small, implying  {$c_{ij} \approx 1$} and $s_{ij}\approx 0$ when the corresponding line is switched on. Therefore, one should try to improve these bounds. 

We  adapt the procedure proposed in \cite{kocuk2015}  (which dealt only with OPF) to the case of OTS in order to obtain variable bounds, that is, we solve a reduced version of the full MISOCP relaxation to efficiently compute bounds. %relaxation.
In particular, to find variable bounds for $c_{kl}$ and $s_{kl}$ for some $(k,l) \in \mathcal{L}$, consider the buses which can be reached from either $k$ or $l$ in at most $r$ steps. Denote this set of buses as $ \mathcal{B}_{kl} (r)$. 
For instance,  $ \mathcal{B}_{kl} (0) = \{k,l\}$, $ \mathcal{B}_{kl} (1) = \delta(k) \cup \delta(l)$, etc. %Similarly, 
{We also define} 
$\mathcal{G}_{kl} (r)=\mathcal{B}_{kl} (r) \cap \mathcal{G}$ and 
$\mathcal{L}_{kl} (r) = 
\{ (i,j) \in \mathcal{L} : i \in \mathcal{B}_{kl} (r) \text{ or }  j \in \mathcal{B}_{kl} (r) \}$. Then, we consider the following SOCP relaxation:
\begin{subequations}\label{bounding SOCP}
\begin{align}
   & p_i^g-p_i^d = g_{ii}c_{ii} + \sum_{j \in \delta(i)} p_{ij}   & i& \in\mathcal{B}_{kl} (r)\\
  &  q_i^g-q_i^d = -b_{ii}c_{ii} + \sum_{j \in \delta(i)}q_{ij}  & i& \in \mathcal{B}_{kl} (r) \\
 & p_{ij} = -G_{ij}c_{ii}^j + G_{ij}c_{ij} -B_{ij}s_{ij}   &(&i,j) \in  \mathcal{L}_{kl} (r) \\
  &q_{ij} = \ \ \,B_{ij}c_{ii}^j -  B_{ij}c_{ij} - G_{ij}s_{ij} &(&i,j) \in  \mathcal{L}_{kl} (r)\\
  & p_{ij}^2+q_{ij}^2  \le (S_{ij}^{\text{max}})^2  &(&i,j) \in \mathcal{L}_{kl} (r)  \\
  & \underline V_i^2 \le c_{ii} \le \overline V_i^2    & i& \in \mathcal{B}_{kl} (r+1) \\
  &   p_i^{\text{min}}  \le p_i^g \le p_i^{\text{max}}     & i& \in \mathcal{G}_{kl} (r) \\
  &   q_i^{\text{min}}  \le q_i^g \le q_i^{\text{max}}     & i& \in \mathcal{G}_{kl} (r) \\
 & \underline c_{ij}x_{ij}\le c_{ij} \le \overline c_{ij} x_{ij}   &(&i,j) \in \mathcal{L}_{kl} (r) \\
  &  \underline s_{ij}x_{ij}\le s_{ij} \le \overline s_{ij} x_{ij}   &(&i,j) \in \mathcal{L}_{kl} (r)\\
  &  \underline c_{ii}x_{ij} \le c_{ii}^j \le \overline c_{ii} x_{ij}   &(&i,j) \in \mathcal{L}_{kl} (r) \\
  &  c_{ii}-\overline c_{ii}(1-x_{ij}) \le c_{ii}^j    &(&i,j) \in \mathcal{L}_{kl} (r) \\
  &   c_{ii}^j \le c_{ii} -  \underline c_{ii}(1- x_{ij})   &(&i,j) \in \mathcal{L}_{kl} (r) \\
  & c_{ij}=c_{ji}, \ \ s_{ij}=-s_{ji}    &(&i,j) \in \mathcal{L}_{kl} (r) \\
  &    c_{ij}^2+s_{ij}^2  \le c_{ii}^jc_{jj}^i     &(&i,j) \in \mathcal{L}_{kl} (r) \\
 & 0 \le x_{ij} \le 1 &(&i,j) \in \mathcal{L}_{kl} (r)\\
 & x_{kl} = 1.  \label{bound socp fixed}
\end{align}
\end{subequations}
{ 
Essentially, \eqref{bounding SOCP} is the continuous relaxation of  MISOCP relaxation applied to the part of the power network within $r$ steps of the buses $k$ and $l$. $c_{kl}$ and $s_{kl}$ can be minimized and maximized subject to \eqref{bounding SOCP} for each edge $(k,l)$ to obtain lower and upper bounds, respectively. These SOCPs can be solved in parallel, since they are independent of each other. It is {observed} that a good tradeoff between accuracy and speed is to select $r=2$\cite{kocuk2015}.  Constraint \eqref{bound socp fixed} may seem to restrict the feasible region, however, the way we defined $c_{kl}$ and $s_{kl}$ variables, they are the values for cosine and sine components when $x_{kl}=1$ (otherwise, they are 0). Therefore, it is enough for the bounds to be valid for $x_{kl}=1$ only.  
}
%In our experiments, larger values of $r$ improve variable bounds marginally.

%{\bf ??? TRIVIAL VALID INEQUALITIES ??? }

Bounds on an artificial edge $(i,j)$ used in the construction of McCormick envelopes are chosen as follows:
\begin{equation}
 \overline c_{ij} = -\underline c_{ij}=  \overline s_{ij} = -\underline s_{ij}= \overline V_i \overline V_j.
\end{equation}

A similar idea can be used to fix some of the binary variables as well. In particular, we can minimize $x_{kl}$ over \eqref{bounding SOCP} after omitting \eqref{bound socp fixed}. If the optimal value turns out to be strictly positive, then $x_{kl}$ can be fixed  to one.

%{\bf Trivial Valid Inequalities}:
%\begin{subequations}
%\begin{align}
%& -S_{ij}^{\text{max}} x_{ij} \le p_{ij}, q_{ij} \le S_{ij}^{\text{max}}x_{ij} \\
%& p_{ij}^2+q_{ij}^2  \le (S_{ij}^{\text{max}} x_{ij})^2 \\
%& c_{ij}^2+s_{ij}^2  \le c_{ii}c_{jj}
%\end{align}
%\end{subequations}

\section{Algorithm}
\label{sec:alg} 

In this section, we propose an algorithm to solve OTS. The algorithm has two {phases}. The first {phase} involves solving a sequence of SOCPs obtained by relaxing integrality restriction of the binary variables in MISOCP \eqref{misocp}, and incorporates cycle inequalities generated from the extended SDP and McCormick relaxations in Section \ref{SOCP SDP sep x} and  \ref{McCormick sep x}. In this {phase}, the aim is to strengthen the lower bound on the MISOCP relaxation. The second {phase} involves solving a sequence of MISOCP relaxations strengthened by cycle inequalities. The aim in this {phase} is to obtain high quality feasible solutions for OTS. In particular, this is achieved by solving OPF subproblems with fixed {topologies obtained from the integral solutions found during the branch-and-cut process of solving the MISOCP \eqref{misocp}}. This procedure is repeated by ``forbidding'' the topologies already considered in order to obtain different network configurations in the subsequent iterations.

Now we formally define  the ingredients of the algorithm. First, let $SOCP(\mathcal{V})$ be the continuous relaxation of MISOCP \eqref{misocp} with a set of valid inequalities $\mathcal{V}$ obtained from cycle inequalities using extended SDP and McCormick relaxations. {The set $\mathcal{V}$ is dynamically updated $T_1$ times}. 
Similarly, we define $MISOCP(\mathcal{V}, \mathcal{F})$ as the MISOCP relaxation of OTS with a set of valid cycle inequalities $\mathcal{V}$ and forbidden topologies $\mathcal{F}$. Here, we forbid a topology $x^*\in\mathcal{F}$ by adding the following ``no-good" cut  (see  \cite{Angulo} for generalizations) to the formulation: 
\begin{equation} \label{no good cut}
\sum_{(i,j):x_{ij}^*=1} (1-x_{ij})  + \sum_{(i,j):x_{ij}^*=0} x_{ij} \ge 1.
\end{equation}
We denote by $LB_t$ as the optimal {value} of $MISOCP(\mathcal{V}, \mathcal{F})$ and $\mathcal{P}_t$ as the set of all integral solutions found by the solver at the $t$-th iteration. For instance, CPLEX offers this option called \textit{solution pool}. {In a practical implementation}, this part is repeated $T_2$ times.

Let $OPF(x)$ denote the value of a feasible solution to OPF problem \eqref{SOCP} for the fixed topology induced by the integral vector $x$. Finally, $UB$ is the best upper bound on OTS. Now, we present Algorithm \ref{alg:acots}.

\begin{algorithm}
\caption{OTS algorithm.}
\label{alg:acots}
\begin{algorithmic}
\STATE Input: $T_1$, $T_2$, $\epsilon$.
\STATE {Phase I}: Set  $\mathcal{V}  \gets \emptyset$, $\mathcal{F} \gets \emptyset$, $UB \gets \infty$.
\FOR{$\tau=1, \dots, T_1$}
\STATE Solve $SOCP(\mathcal{V})$.
\STATE Separate cycle inequalities for each cycle in a cycle basis to obtain a set of valid inequalities $\mathcal{V}_\tau$.
\STATE Update $\mathcal{V} \gets \mathcal{V}\cup \mathcal{V}_{\tau}$.
\ENDFOR
\STATE {Phase II}: Set {$t\gets 0$.}
\REPEAT 
\STATE $t \gets t+1$
\STATE Solve $MISOCP(\mathcal{V}, \mathcal{F})$ to obtain a pool of integral solutions $\mathcal{P}_t$ and record the optimal {cost} as $LB_t$.
\FORALL{$x \in \mathcal{P}_t$}
\IF{ $OPF(x) < UB$ } 
	\STATE $UB \gets OPF(x)$
\ENDIF
\ENDFOR
\STATE Update $\mathcal{F} \gets \mathcal{F}  \cup \mathcal{P}_t$.
\UNTIL{{${LB}_t \ge (1 - \epsilon)  UB$  or  $t \ge T_2$}}
\end{algorithmic}
\end{algorithm}

\begin{obs}\label{ideal alg}
If $OPF(x)$ returns globally optimal solution for every topology $x$, $\epsilon=0$ and $T_2 = \infty$, then Algorithm~\ref{alg:acots} converges to the optimal solution of  OTS in finitely many iterations.
\end{obs}

{Observation}  \ref{ideal alg} follows from the fact that there are finitely many topologies and by the hypothesis that $OPF(x)$ can be solved globally{, which is possible for some IEEE instances using moment/sum-of-squares relaxations \cite{Josz}}. Although {Observation}   \ref{ideal alg} states that Algorithm \ref{alg:acots} can be used to solve OTS to global optimality in finitely many iterations, the requirement of solving  $OPF(x)$ to global optimality may not be satisfied {always}. In practice, we can  solve OPF subproblems using local solver, in which case we have {Observation} \ref{practical alg}.

\begin{obs}\label{practical alg}
If $OPF(x)$ is solved by  a local solution method, then we have ${{LB}_1} \le  z^* \le UB$ upon termination of Algorithm \ref{alg:acots}, where $z^*$ is the optimal value of OTS.
\end{obs}

\section{Computational Experiments}
\label{sec:comp expr}

In this section, we present the results of our extensive computational experiments on standard IEEE instances available from MATPOWER \cite{Matpower} and instances from NESTA 0.3.0 archive with congested operating conditions \cite{nesta}. The code is written in the C\# language with Visual Studio 2010 as the compiler. For all experiments, we used a 64-bit computer with Intel Core i5 CPU 2.50GHz processor and 16 GB RAM. Time is measured in seconds. 
We use three different solvers:
\begin{itemize}
\item
CPLEX 12.6 \cite{Cplex126}  to solve MISOCPs.
\item
Conic interior point solver MOSEK 7.1 \cite{MOSEK} to solve  SDP separation problems.
\item
Nonlinear interior point solver IPOPT \cite{wachter} to find local optimal solutions to $OPF(x)$.
\end{itemize}
We use a Gaussian elimination based approach to construct a cycle basis proposed in \cite{kocuk2014Switch} and use this set of cycles  in the separation phase.

\subsection{Methods}
%In our preliminary examination, we have tried several different SOCP relaxations and decided to report the results of the following four settings:
We report the results of three algorithmic settings:
\begin{itemize}
\item
$\mathsf{SOCP}$:  MISOCP formulation \eqref{misocp} {in Phase II} without {Phase I (i.e. $T_1=0$)}. %our base  SOCP  without any improvements.
\item
{$\mathsf{SOCPA}$}: $\mathsf{SOCP}$ strengthened {by the} arctangent envelopes {introduced in Section \ref{sec:arctan envelope x}}.
%\item
%{$\mathsf{SOCPA\_SDP}$}: $\mathsf{SOCPA}$ strengthened further {by dynamically generating linear} valid inequalities obtained from separating an SOCP feasible solution from the SDP relaxation over cycles.  conv$(\mathcal{S}_1 \cup \mathcal{S}_0)$
%%The separation routine is developed in Section \ref{SOCP SDP sep}.} %submatrices corresponding to cycles by solving SDP separation problem \eqref{sdp sep}.
\item
{$\mathsf{SOCPA\_Disj}$}: $\mathsf{SOCPA}$ strengthened further {by dynamically generating linear} valid inequalities obtained from separating an SOCP feasible solution from the SDP and McCormick relaxation over cycles using a disjunctive argument $T_1$ times. In particular, a separation oracle is used to separate a given
point from  conv$( (\mathcal{S}_1 \cap \mathcal{M}_1) \cup \mathcal{S}_0)$.
\end{itemize}

The following four performance measures are used to assess the accuracy and the efficiency of the proposed methods:
\begin{itemize}
\item
``\%OG" is the percentage optimality gap proven by our algorithm calculated as $100 \times (1-\widetilde{LB}_1 /UB)$. Here,  $\widetilde{LB}_1$ is the lower bound proven, which may be strictly smaller than ${LB}_1$ due to optimality gap tolerance and time limit.
\item
``\%CB" is the percentage cost benefit obtained by line switching calculated as $100 \times (1-UB / OPF(e))$, where $e$ is the vector of ones so that $OPF(e)$ corresponds to the OPF solution with the initial topology.
\item
``\#off" is the number of lines switched off in the topology which gives $UB$.
\item
``TT" is the total time in seconds, {including} preprocessing (bound tightening), solution of $T_1=5$ {rounds of} SOCPs to improve lower bound and separation problems to generate cutting planes (in the case of {$\mathsf{SOCPA\_Disj}$}), solution of $T_2$ {rounds of} MISOCPs and several calls to local solver IPOPT with given topologies.   MISOCPs   are solved under 720 seconds time limit so that 5 iterations take about 1 hour (optimality gap for integer programs is 0.01\%).  Preprocessing and separation subproblems are parallelized.
\end{itemize}
We choose parameter $T_2=5$ and pre-terminate Algorithm \ref{alg:acots} if 0.1\% optimality gap is proven.

\subsection{Results}

The results of our computational experiments are presented in Tables \ref{table:IEEE} and \ref{table:nesta} for standard IEEE and NESTA instances, respectively.   We considered instances up to 300-bus since Phase II of the Algorithm \ref{alg:acots} does not scale up well for larger instances. Let us start with the former: IEEE instances are a relatively easy set since transmission line limits are generally not binding. Therefore, cost benefits obtained by switching are also limited. The largest cost reduction is obtained for case30Q with 2.24\%.  Among the three methods, the most successful 
one  is $\mathsf{SOCPA\_Disj}$, {on average} proving 0.05\% optimality gap and providing 0.31\% cost savings. In terms of computational {time}, $\mathsf{SOCP}$ is the fastest, however, its performance is not as good as the other two. Quite interestingly, $\mathsf{SOCPA\_Disj}$ is faster than $\mathsf{SOCPA}$, on average, for this set of instances. % \textcolor{blue}{mainly due to the longer computation time on 29edin and 300ieee instances.}
%%all cases solved up to 0.17\% of global optimal with $\mathsf{SOCPA\_Disj}$ with an average of 0.05\%
{In terms of comparison with other methods,} unfortunately, there is limited literature {for this purpose}. In \cite{hijazi2013}, nine of these instances (except for cases 9Q and 30Q) are considered and a quadratic convex (QC) relaxation based approach is used. On average, their approach proves 0.14\% optimality gap, which is worse than any of our methods over the same nine instances. The only instance QC approach is better is {118ieee} with 0.11\% optimality gap, while it is worse than our methods for case300 with a 0.47\% optimality gap.

%\addtolength{\tabcolsep}{-.15pt}    
% Table generated by Excel2LaTeX from sheet 'Tables for Paper'
\begin{table}[H]\footnotesize
\begin{center}
\caption{Results summary for standard IEEE Instances.}
\label{table:IEEE}
% Table generated by Excel2LaTeX from sheet 'Tables for Paper'
\begin{tabular}{c|rrrr|rrrr|rrrr}

       \cline{2-13}
    &                         \multicolumn{ 4}{c|}{SOCP} &                        \multicolumn{ 4}{c|}{SOCPA} &                \multicolumn{ 4}{c}{SOCPA\_Disj} \\

\hline
      case &       \%OG &       \%CB &       \#off &         TT(s) &       \%OG &       \%CB &      \#off &         TT(s) &      \%OG &       \%CB &       \#off &         TT(s) \\
\hline
       6ww &       0.16 &       0.48 &          2 &       1.29 &       0.02 &       0.48 &          2 &       0.67 &        0.01 &       0.48 &          2 &       1.28 \\

         9 &       0.00 &       0.00 &          0 &       0.26 &       0.00 &       0.00 &          0 &       0.22 &           0.00 &       0.00 &          0 &       0.55 \\

        9Q &       0.04 &       0.00 &          0 &       0.42 &       0.04 &       0.00 &          0 &       0.33 &            0.04 &       0.00 &          0 &       0.97 \\

        14 &       0.08 &       0.00 &          0 &       0.66 &       0.09 &       0.00 &          0 &       0.70 &           0.01 &       0.00 &          1 &       1.81 \\

    ieee30 &       0.05 &       0.00 &          1 &       1.95 &       0.05 &       0.00 &          0 &       1.67 &          0.02 &       0.00 &          1 &       3.84 \\

        30 &       0.07 &       0.52 &          1 &       4.60 &       0.06 &       0.52 &          2 &       5.01 &           0.03 &       0.51 &          2 &       9.39 \\

       30Q &       0.44 &       2.05 &          2 &      24.43 &       0.43 &       2.03 &          5 &      25.80 &        0.13 &       2.24 &          5 &      44.16 \\

        39 &       0.03 &       0.00 &          0 &       2.53 &       0.01 &       0.02 &          1 &       3.17 &        0.01 &       0.02 &          1 &       4.48 \\

        57 &       0.07 &       0.02 &          4 &       6.18 &       0.07 &       0.02 &          4 &       8.72 &        0.08 &       0.01 &          1 &      13.59 \\

       118 &       0.19 &       0.08 &          4 &    3065.64 &       0.15 &       0.12 &         10 &    2553.59 &           0.17 &       0.08 &         16 &    3174.01 \\

       300 &       0.16 &       0.02 &          9 &    2318.89 &       0.15 &       0.03 &         12 &    3624.12 &            0.10 &       0.05 &         15 &    2803.31 \\

   {\bf  Avg.} & {\bf 0.12} & {\bf 0.29} & {\bf 2.1} & {\bf 493.35} & {\bf 0.10} & {\bf 0.29} & {\bf 3.3} & {\bf 565.82}  & {\bf 0.05} & {\bf 0.31} & {\bf 4.0} & {\bf 550.67} \\
\hline
\end{tabular}  
\end{center}
\end{table}
%\addtolength{\tabcolsep}{.15pt}    

%\addtolength{\tabcolsep}{-.5pt}    
\begin{table}[H]\footnotesize
\begin{center}
\caption{Results summary for NESTA Instances from Congested Operating Conditions.}
\label{table:nesta}
% Table generated by Excel2LaTeX from sheet 'Tables for Paper'
\begin{tabular}{c|rrrr|rrrr|rrrr}

       \cline{2-13}
    &                         \multicolumn{ 4}{c|}{SOCP} &                        \multicolumn{ 4}{c|}{SOCPA}  &               \multicolumn{ 4}{c}{SOCPA\_Disj} \\
\hline
      case &       \%OG &       \%CB &       \#off &         TT(s) &           \%OG &       \%CB &       \#off &         TT(s) &       \%OG &       \%CB &       \#off &         TT(s) \\
\hline
     3lmbd &       3.30 &       0.00 &          0 &       0.14 &       2.00 &       0.00 &          0 &       0.14 &            1.17 &       0.00 &          0 &       0.30 \\

       4gs &       0.65 &       0.00 &          0 &       0.11 &       0.16 &       0.00 &          0 &       0.13 &            0.00 &       0.00 &          0 &       0.27 \\

      5pjm &       0.18 &       0.27 &          1 &       0.61 &       0.01 &       0.27 &          1 &       0.41 &           0.02 &       0.27 &          1 &       0.89 \\

       6ww &       6.06 &       7.74 &          1 &       1.23 &       1.34 &       7.74 &          1 &       1.64 &          1.05 &       7.74 &          1 &       1.97 \\

     9wscc &       0.00 &       0.00 &          0 &       0.19 &       0.00 &       0.00 &          0 &       0.20 &           0.00 &       0.00 &          0 &       0.30 \\

    14ieee &       1.02 &       0.33 &          1 &       2.86 &       0.89 &       0.45 &          2 &       3.48 &           0.41 &       0.45 &          2 &       4.49 \\

    29edin &       0.43 &       0.00 &          2 &      12.79 &       0.24 &       0.18 &         13 &     299.82 &          0.33 &       0.08 &         21 &     181.74 \\

      30as &       1.81 &       3.13 &          2 &      14.82 &       0.35 &       3.30 &          5 &      19.52 &         0.34 &       3.30 &          5 &      24.93 \\

     30fsr &       3.24 &      44.20 &          2 &       9.72 &       0.05 &      44.98 &          2 &       4.76 &            0.03 &      44.98 &          3 &       6.97 \\

    30ieee &       0.54 &       0.46 &          1 &      12.28 &       0.40 &       0.48 &          2 &      10.61 &             0.15 &       0.48 &          2 &      13.37 \\

    39epri &       1.92 &       1.10 &          1 &      11.56 &       0.80 &       1.41 &          2 &      13.20 &            0.70 &       1.52 &          2 &      12.65 \\

    57ieee &       0.12 &       0.10 &          3 &      41.48 &       0.12 &       0.10 &          2 &      58.97 &           0.09 &       0.10 &          3 &      29.86 \\

   118ieee &      41.67 &       4.33 &          3 &     225.57 &      21.51 &      27.98 &         30 &    3838.62 &            7.50 &      39.09 &         21 &    3856.76 \\

   162ieee &       0.57 &       1.05 &          9 &    3675.75 &       0.63 &       1.00 &         15 &    3861.29 &           0.60 &       1.00 &         15 &    3855.50 \\

   189edin &       5.31 &       1.10 &          3 &     540.02 &       4.81 &       0.13 &          2 &    2194.80 &         5.58 &       0.00 &          0 &    3634.02 \\

   300ieee &       1.00 &       0.10 &         12 &    3655.10 &       0.65 &       0.37 &         21 &    3640.14 &          0.61 &       0.35 &         21 &    3651.95 \\

 {\bf  Avg.} & {\bf 4.24} & {\bf 3.99} & {\bf 2.6} & {\bf 512.76} & {\bf 2.12} & {\bf 5.52} & {\bf 6.1} & {\bf 871.73}  & {\bf 1.16} & {\bf 6.21} & {\bf 6.1} & {\bf 954.75} \\
\hline
\end{tabular}  

\end{center}
\end{table}
%\addtolength{\tabcolsep}{.5pt}    

Now let us consider NESTA instances with congested operating conditions. This set is particularly suited for line switching as more stringent transmission line limits are imposed. In fact, large cost improvements are observed for some test cases. For instance, about 45\% and  39\% cost reductions are possible for cases 30fsr and 118ieee, respectively. Other instances with sizable cost reductions include cases 6ww and 30as. $\mathsf{SOCPA\_Disj}$ is again the most successful method if we look at averages of optimality gap (1.16\%) and cost savings (6.21\%). It certifies that the best topology is within 1.17\% of the optimal for all the cases except for 118ieee and 189edin. In terms of computational {time}, $\mathsf{SOCP}$ is again the fastest, however, its performance is significantly worse than the other two. We also note that  $\mathsf{SOCPA}$ improves quite a bit over  $\mathsf{SOCP}$ in terms of optimality and cost benefits with 70\% increase in computational time.  $\mathsf{SOCPA\_Disj}$ takes about only 10\% more time than  $\mathsf{SOCPA}$. As we go from $\mathsf{SOCP}$ to $\mathsf{SOCPA\_Disj}$, problems get more complicated and sometimes, MISOCPs are not solved to optimality within time limit. Consequently, for cases 189edin and 300ieee, the optimality gaps proven and cost benefits obtained by  $\mathsf{SOCPA\_Disj}$ can be slightly worse.

Finally, we note that that optimality gaps can be explained by two non-convexities:
1) integrality, 2) power flow equations. For instance, in case 3lmbd, the {optimality} gap can only be explained by the non-convexity of power flow equation since all the relevant topologies are considered. Similarly, at least some portion of the relatively large optimality gaps for cases 118ieee and 189edin may be attributed to  non-convexity of power flow equations. Consequently, any future improvements on strengthening the convex relaxations of OPF problem can  be useful in closing more gaps in OTS as well.

{
\subsection{Discussion}
In this section, we take a closer look at some of the instances with large cost benefits and try to gain some insight as to 1) how these large savings are obtained, and 2) how simple heuristics may fail to produce comparable results. 
 Firstly, using a small example, we illustrate how large cost savings can be obtained. 
Secondly, we  compare the results of our algorithm with a commonly used heuristic based on switching the best line and demonstrate how different the  solution quality can be.

To address the first issue, let us concentrate on a small instance, namely case6ww from NESTA archive. This instance has the same topology and line characteristics as the standard IEEE test case but load and generation parameters are slightly different. In particular, $p_i^{d}=78.24$, $q_i^{d}=70$, $\underline V_i = 0.95$ and $\overline V_i=1.05$ for the load buses $i=4,5,6$ while the data for generation buses 1, 2 and 3 is summarized in Table \ref{GenData6ww}.
\begin{table}[h!]
\caption{Generator data for NESTA case6ww test case.}
\label{GenData6ww}
\begin{center}
\begin{tabular}{rrrrrrr}
\hline
 &          $p_i^{\text{min}}$ & $p_i^{\text{max}}$     &      $q_i^{\text{min}}$ & $q_i^{\text{max}}$ & $\underline V_i =\overline V_i $   &       cost \\
\hline
         $1$ &               $25$ &        $200$ &        $-100$ &        $100$ &     $ 1.05$&    $1.276311$ \\    
         $2$ &               $18.75$ &        $106$ &        $-100$ &        $100$ &   $1.05$  &   $0.586272$ \\
         $3$ &               $22.5$ &        $93$ &        $-100$ &        $100$ &   $1.07$ &   $1.29111$ \\
\hline
\end{tabular}
\end{center}
\end{table}
With this topology, the local optimal solution  obtained using IPOPT  with objective value of 273.76 is given in Figure \ref{case6wwSOL}. We note that the lines $(1,5)$, $(2,4)$ and $(3,6)$ are congested in this configuration. On the other hand, if the line $(1, 2)$ is switched off, then the objective value reduces to 252.57, corresponding to a $7.74\%$ cost saving over the initial topology. The difference is due to the fact that the outputs of generators 1 and 2 are now changed to  (85.56, 32.74) and  (84.25, 63.26), respectively.  Notice that with  the new topology, the cheaper generator 2 is used more, which results in the cost reduction. In the initial topology, this is not possible since the voltage magnitudes of the generators are fixed, and lines $(1,5)$ and $(2,4)$ are congested.

{

%%https://github.com/rwl/thesis/blob/master/tikz/case6ww.tex
%http://academia.stackexchange.com/questions/14010/how-do-you-cite-a-github-repository

\newcommand{\genset}[3]{
  \node[circle,draw,thick,minimum width=7mm,inner sep=0pt,fill=white,#3] (#1) at (#2) {};
  \draw[thick] ($(#2)-(2mm,0)$) sin ++(1mm,1mm) cos ++(1mm,-1mm) sin
  ++(1mm,-1mm) cos ++(1mm,1mm); }

\begin{figure}[H]
\centering
\caption{
{Flow diagram for the solution of NESTA case6ww without any line switching. The numbers above each generator node respectively represent the active and reactive power  output. Similarly, the  numbers  near each edge respectively represent the active and reactive power flow of the line in the direction from the small indexed bus to the large indexed one. The figure is generated by modifying \cite{Lincoln}.}
} 
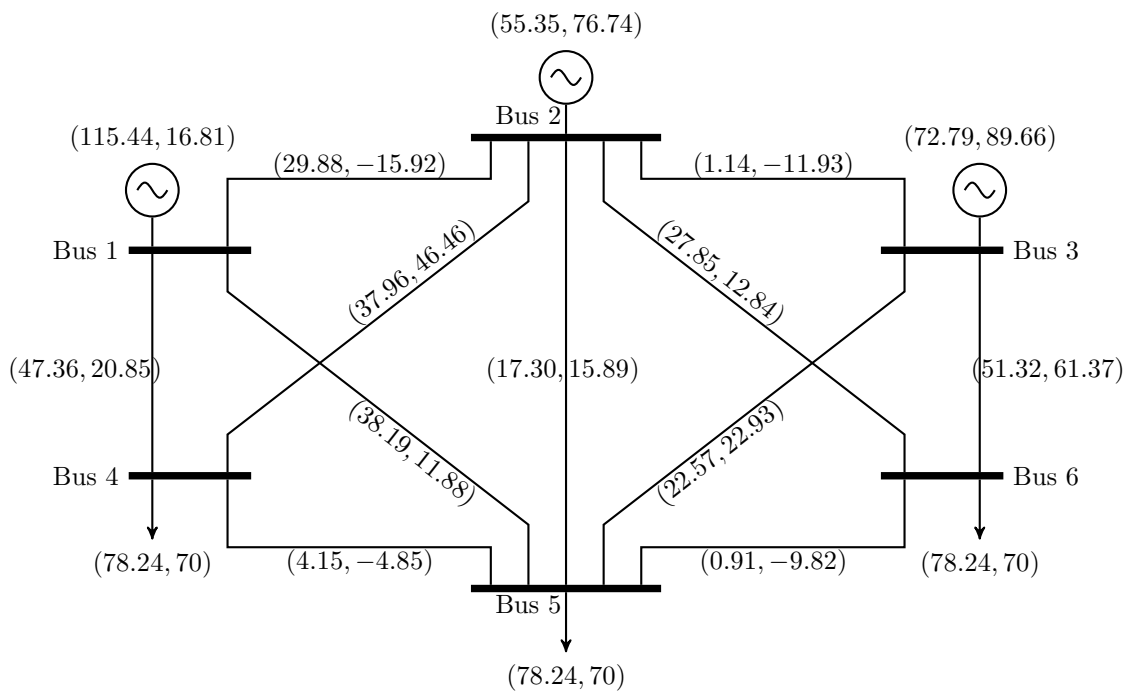
\label{case6wwSOL}
\begin{tikzpicture}[thick]

\tikzstyle{line}=[-, thick]
\tikzstyle{loadline}=[->,thick,>=stealth']
\tikzstyle{busbar} = [rectangle,draw,fill=black,inner sep=0pt];
\tikzstyle{hbus} = [busbar,minimum width=10mm,minimum height=2pt];

\coordinate (c1) at (0,3);
\coordinate (c2) at (5,4.5);
\coordinate (c3) at (10,3);
\coordinate (c4) at (0,0);
\coordinate (c5) at (5,-1.5);
\coordinate (c6) at (10,0);
\coordinate (over) at (0,3.75);

\begin{footnotesize}
\node[hbus,minimum width=16mm,label=left:Bus 1] (b1) at (c1) {};
\node[hbus,minimum width=25mm,label={[xshift=-0.5cm, yshift=0.0cm]Bus 2}] (b2) at (c2) {};
\node[hbus,minimum width=16mm,label=right:Bus 3] (b3) at (c3) {};
\node[hbus,minimum width=16mm,label=left:Bus 4] (b4) at (c4) {};
\node[hbus,minimum width=25mm,label={[xshift=-0.5cm, yshift=-0.5cm]Bus 5}] (b5) at (c5) {};
\node[hbus,minimum width=16mm,label=right:Bus 6] (b6) at (c6) {};
\end{footnotesize}

% \busbar{b1}{c1}{20mm}
% \busbar{b2}{c2}{40mm}
% \busbar{b3}{c3}{20mm}
% \busbar{b4}{c4}{20mm}
% \busbar{b5}{c5}{40mm}
% \busbar{b6}{c6}{20mm}

% Branch 1-2.
\draw[line] ([xshift=5mm] b1.north) |- ([xshift=-10mm,yshift=-5mm] b2.south) --([xshift=-10mm] b2.south) 
node[text centered, text width=30mm, xshift=-1.75cm, yshift=-0.325cm] {\footnotesize$(29.88,	-15.92)$ } ;
% Branch 2-3.
\draw[line] ([xshift=10mm] b2.south) 
node[text centered, text width=30mm, xshift=1.75cm, yshift=-0.325cm]{\footnotesize$(1.14,	-11.93)$} -- ++(0,-0.5) -| ([xshift=-5mm] b3.north) ;
% Branch 4-5.
\draw[line] ([xshift=5mm] b4.south) |- ([xshift=-10mm,yshift=5mm] b5.north) --([xshift=-10mm] b5.north)
node[text centered, text width=30mm, xshift=-1.75cm, yshift=0.3cm] {\footnotesize$(4.15,	-4.85)$};
% Branch 5-6.
\draw[line] ([xshift=10mm] b5.north)
node[text centered, text width=30mm, xshift=1.75cm, yshift=0.3cm]{\footnotesize$(0.91,	-9.82)$} -- ([xshift=10mm,yshift=5mm] b5.north) -|([xshift=-5mm] b6.south);

% Branch 1-4.
\draw[line] ([xshift=-5mm] b1.south) -- ([xshift=-5mm] b4.north) node[text centered, text width=30mm, xshift=-0.9cm, yshift=1.35cm] {\footnotesize$(47.36,	20.85)$};
% Branch 2-5.
\draw[line] (b2.south) -- (b5.north) node[text centered, text width=30mm, xshift=-0.05cm, yshift=2.85cm] {\footnotesize$(17.30,	15.89)$};
% Branch 3-6.
\draw[line] ([xshift=5mm] b3.south) -- ([xshift=5mm] b6.north) node[text centered, text width=30mm, xshift=0.9cm, yshift=1.35cm] {\footnotesize$(51.32,	61.37)$};

% Branch 1-5.
\draw[line] ([xshift=5mm] b1.south) -- ([xshift=5mm,yshift=-5mm] b1.south)-- ([xshift=-5mm,yshift=8mm] b5.north) -- ([xshift=-5mm] b5.north) 
node[text centered, text width=30mm, xshift=-1.55cm, yshift=1.75cm,rotate=-38]{\footnotesize$(38.19,	11.88)$};
% Branch 3-5.
\draw[line] ([xshift=-5mm] b3.south) -- ([xshift=-5mm,yshift=-5mm] b3.south)-- ([xshift=5mm,yshift=8mm] b5.north) -- ([xshift=5mm] b5.north)
node[text centered, text width=30mm, xshift=1.55cm, yshift=1.75cm,rotate=38]{\footnotesize$(22.57,	22.93)$};
% Branch 2-4.
\draw[line] ([xshift=-5mm] b2.south) 
node[text centered, text width=30mm, xshift=-1.55cm, yshift=-1.75cm,rotate=38]{\footnotesize$(37.96,	46.46)$}
--([xshift=-5mm,yshift=-8mm] b2.south) -- ([xshift=5mm,yshift=5mm] b4.north)-- ([xshift=5mm] b4.north);
% Branch 2-6.
\draw[line] ([xshift=5mm] b2.south) 
node[text centered, text width=30mm, xshift=1.55cm, yshift=-1.75cm,rotate=-38]{\footnotesize$(27.85,	12.84)$}
-- ([xshift=5mm,yshift=-8mm] b2.south)-- ([xshift=-5mm,yshift=5mm] b6.north) -- ([xshift=-5mm] b6.north);

% Generator 1.
\genset{g1}{$(c1)+(-5mm,8mm)$}{label=above:{\footnotesize{$(115.44,	16.81)$}}}
\draw[line] ([xshift=-5mm] b1.north) -- (g1.south);
% Generator 2.
\genset{g2}{$(c2)+(0,8mm)$}{label=above:{\footnotesize{$(55.35,	76.74)$}}}
\draw[line] (b2.north) -- (g2.south);
% Generator 3.
\genset{g3}{$(c3)+(5mm,8mm)$}{label=above:{\footnotesize{$(72.79,	89.66)$}}}
\draw[line] ([xshift=5mm] b3.north) -- (g3.south);

% Load 1.
\draw[loadline] ([xshift=-5mm] b4.south) -- ++(0,-0.8) node[text centered,text
width=20mm,below] {\footnotesize{$(78.24, 70)$}};
% \loadd{l1}{$(c4)-(5mm,15mm)$}
% \draw[line] (l1.south) -- ([xshift=-5mm] b4.south);
% Load 2.
\draw[loadline] (b5.south) -- ++(0,-0.8) node[text centered,text
width=20mm,below] {\footnotesize{$(78.24, 70)$}};
% \loadd{l2}{$(c5)-(0mm,15mm)$}
% \draw[line] (l2.south) -- (b5.south);
% Load 3.
\draw[loadline] ([xshift=5mm] b6.south) -- ++(0,-0.8) node[text centered,text
width=20mm,below] {\footnotesize{$(78.24, 70)$}};
% \loadd{l3}{$(c6)+(5mm,-15mm)$}
% \draw[line] (l3.south) -- ([xshift=5mm] b6.south);

\end{tikzpicture}
\end{figure}

}

Now, let us consider the second issue. Due to the combinatorial nature of OTS problem, heuristics are frequently used to obtain suboptimal solutions. A commonly used one is to switch off a single line to obtain cost benefits \cite{Sahraei2014, Korad2015}. Although this heuristic idea is easy to implement and works well in some instances, there are no guarantees on its accuracy. For example, in case6ww, the best line to switch off is, in fact, $(1, 2)$ suggested by both the best line heuristic and  our algorithm. However, for other problems with large cost benefits, this is not always the case. For instance, in case30as and case30fsr, the best line heuristic reduces the overall cost to $2.99\%$ and $44.02\%$ respectively, compared to  $3.30\%$ and $44.98\%$ obtained from our algorithm. For case118ieee, the cost reduction is dramatically different. The best line heuristic reduces the cost only by $19.52\%$ while our algorithm provides a topology with 39.09\% saving. Moreover, the best line heuristic does not provide any guarantee on how good the solution is while our algorithm gives optimality guarantees by construction. Therefore, adapting Algorithm \ref{alg:acots} in real-time operations can yield significant savings over simple heuristics.
}

\section{Conclusions}
\label{sec:concl}

%In this paper, we proposed a systematic approach to study the AC OTS problem. 
%%Our work originates from the recent developments for the closely related AC OPF problem. 
%Our work include an alternative formulation for OTS and its MISOCP relaxation. We improve this relaxation via arctangent envelopes and cutting planes obtained using disjunction.  The use of disjunctive cuts help improve gap closure significantly. 
%%To the best of our knowledge, such SDP disjunctive cuts have never been computationally used before. 
%%
%We also propose a practical algorithm to obtain high quality feasible OTS solutions. Our experiments on standard and congested instances suggest that the proposed methods are effective in obtaining strong relaxations and provably good feasible solutions. To the best of our knowledge, the results obtained in this paper outperform the findings on the limited literature available.
%%
%We remind the reader that AC OTS is an especially challenging problem since it embodies two types of nonconvexities due to AC power flows and integral variables. We expect that the methodology developed in this paper will help solve AC OTS problem in real life operations.

In this paper, we proposed a systematic approach to solve the AC OTS problem. In particular, we presented an alternative formulation for OTS and constructed a MISOCP relaxation. We improved the strength of this relaxation by the addition of arctangent envelopes and cutting planes obtained using disjunctive techniques. The use of these disjunctive cuts help in closing gap significantly. Our experiments on standard and congested instances suggest that the proposed methods are effective in obtaining strong lower bounds and producing provably good feasible solutions. 
%To the best of our knowledge, on an average the algorithm presented in this paper out performs all other algorithm presented in the literature. 

We remind the reader that AC OTS is a challenging problem since it embodies two types of non-convexities due to AC power flow constraints and integrality of variables. We hope that the methodology developed in this paper can eventually be further improved to solve AC OTS problem in real life operations.
{
As a future work, we would like to pursue finding ways to improve the solution time of MISOCPs as this step is the bottleneck in Algorithm \ref{alg:acots}. Also, decomposition methods can be sought to solve large-scale problems more efficiently, which could make the proposed approach adaptable to real life instances.
}

\vspace{-0mm}

\appendix

\section{Convex Hull of Union of Two Conic Representable Sets} \label{app:union}
Let ${S}_1$ and ${S}_2$ be two bounded, conic representable sets %given by
\begin{equation*}
{S}_i = \{x : \exists u^i: A_i x  + B_i u^i\succeq_{K_i} b_i \} \quad i = 1,2.
\end{equation*}
Here, $K_i$'s are regular (closed, convex, pointed with non-empty interior) cones. 
Then, a conic representation for conv$({S}_1 \cup {S}_2)$ is given as follows:
\begin{subequations}
\begin{align*}
& x = x^1 + x^2,  \, \lambda_1+\lambda_2 = 1, \, \lambda_1, \lambda_2 \geq 0 \\
&  A_i x^i  + B_i u^i\succeq_{K_i} b_i \lambda_i \quad i = 1,2.
\end{align*}
\end{subequations}

\vspace{-0mm}

\section{Separation from an Extended Conic Representable Set} \label{app:sep}
Let ${S}$ be a conic representable set ${S} = \{x : \exists u: Ax  + B u\succeq_{K} b \}$. %given by
%\begin{equation}
%{S} = \{x : \exists u: Ax  + B u\succeq_{K} b \}.
%\end{equation}
Here, $K$ is a regular cone. {Suppose we want to decide if a given point $x^*$ belongs to ${S}$ and find a separating hyperplane $\alpha^\top x\ge\beta$ if $x^*\notin S$. This problem can be formulated as $\max_{\alpha,\beta}\left\{\beta-\alpha^\top x^* : \alpha^\top x\geq \beta \; \forall x\in S \right\}$, where the constraint can be further dualized as  
\begin{align*}
Z^*:=\max_{\alpha,\beta,\mu}\{ \beta-\alpha^\top x^* \, : \; & b^\top\mu\geq\beta, A^\top\mu=\alpha, B^\top\mu=0,\\
                       & \mu\in K^*, -e \le \alpha \le e, -1 \le  \beta \le 1 \},
\end{align*}
where $K^*$ is the dual cone of $K$. If $Z^*\leq 0$, then $x^*\in S$, otherwise, the optimal $\alpha,\beta$ from the above program gives the desired separating hyperplane. For details, please see \cite{kocuk2015}.}
%Let us consider the following  problem
%\begin{subequations}\label{eq:sep cone all}
%\begin{align}
%v^*:=  \max_{\mu, \alpha}  &\hspace{0.5em}  \beta - \alpha^T x^*\\
%\mathrm{s.t.}   &\hspace{0.5em} A^T \mu = \alpha \label{eq:sep cone1} \\
%&\hspace{0.5em} B^T \mu = 0 \\
%&\hspace{0.5em} b^T \mu \ge \alpha \\
%&\hspace{0.5em} \mu \in K^*  \label{eq:sep cone2} \\
%&\hspace{0.5em} -e \le \alpha \le e, -1 \le  \beta \le 1, \label{eq:sep cone bound}
%\end{align}
%\end{subequations}
%where $K^*$ is the dual cone. Here, \eqref{eq:sep cone1}-\eqref{eq:sep cone2} is the dual system equivalent to the condition that $\alpha^T x\leq \beta$ for all $x\in S$; \eqref{eq:sep cone bound} bounds the coefficients $\alpha, \beta$. If $v^*>0$, then the corresponding optimal solution $(\alpha,\beta)$ of \eqref{eq:sep cone all} gives a separating hyperplane such that $\alpha^T x^*<\beta$ and $\alpha^T x \ge \beta$ for all $x \in S$. If $v^*\le 0$, then \eqref{eq:sep cone all} certifies that $x$ is contained in $S$ for some $u$. 

\vspace{-0mm}

\bibliographystyle{plain}
\bibliography{references}

\end{document}